\documentclass{article}
\usepackage{graphicx}
\usepackage{natbib}

\usepackage{amsmath}  
\usepackage{wrapfig}	    
\usepackage{tikz}
\usetikzlibrary{calc}

\begin{document}

\title{On implementation of the minimal perimeter triangle enclosing a convex polygon algorithm
}


\author{Victor~Ermolaev
}




\maketitle

\begin{abstract}
A high-level description of an algorithm which computes the minimum perimeter triangle enclosing
a convex polygon in linear time exists in the literature. Besides that an implementation of the algorithm
is given in the subsequent work. However, that implementation is incomplete and omits
certain common cases. This note's contribution consist of a modern treatment of the algorithm
and its consequent simple and robust implementation.
\end{abstract}

\section{Introduction}
A very simple and ingenious algorithm of finding a minimal perimeter triangle enclosing
a convex polygon was given by~\citet{bhattacharya2002minimum}. The algorithm
has a linear complexity w.r.t. the number of edges of the polygon. This makes it very attractive for applications.

Let a convex polygon be given, such that no two neighbouring edges are colinear.
The algorithm suggests to treat every edge sequentially as one of the sides of the optimal
triangle; the second (suboptimal) side is bootstrapped and one optimization step is performed.
For the exposition of the algorithm suppose that two sides of the minimal perimeter enclosing triangle
are found and form a wedge. A quest for the third side is an optimization procedure w.r.t. edges and vertexes of the polygon:
the third side will either be colinear to one of the edges or pass through
one of the vertexes. Findings  in \citep{bhattacharya2002minimum} state that:
\begin{itemize}
\item if the wedge and extension of a polygon edge form a triangle, then the third side
will be colinear to the edge if a circle escribed to the triangle touches the edge
(and not its extension);
\item if a circle is inscribed into the wedge such that that it passes through one
of the polygon vertexes, then the third side will be a tangent
to the circle at the vertex provided it is also tangent to the polygon.
\end{itemize}
Finding such a side embodies one optimization step of ``closing'' the wedge in the optimal way.
The next optimization step is performed for a wedge formed by the first side and
the newly found side. Such wedge ``flipping'' continues until no further improvement
in perimeter is possible.

As described, the algorithm relies on four subsidiary procedures of fitting a circle into
a wedge while satisfying some constraints. The first two sub-problems allow us to bootstrap
the algorithm, whereas the last two constitute one optimization step. For the bootstrapping part,
call a wedge degenerate, if its arms are parallel.
\paragraph{\textbf{Subsidiary problems}}
\begin{enumerate}
	\item Given a degenerate wedge and a line crossing it, find a circle touching both
		of the wedge's arms and the line;
	\item given a degenerate wedge and a point between its arms, find
		a circle touching the arms of the wedge and passing through the point;
	\item given a wedge and a line crossing it, find an excircle for
		the generated triangle such it touches the given line;
	\item given a wedge and a point between its arms, find
		a circle inscribed into the wedge and passing through the point.
\end{enumerate}
\citet{medvedeva2003implementation} sketch an implementation which is
further developed in~\citep{medvedeva2003computing}. That implementation relies on
slope-intercept representation of lines, which leads to the numerous case distinction.
We develop an alternative treatment of the sub-problems employing simple linear algebra
techniques. Our approach eliminates the infinite-slope problem while
leaving the main body of the algorithm intact. Our proposition reduces the overall
complexity of the algorithm leading to a simpler and more direct implementation,
for example, on top of a linear algebra library. In particular, we provide a reference to 
our implementation hosted on GitHub~\citep{github-release}. 
The GitHub project includes sample polygons from~\citep{parvu2014implementation} as test cases.

\section{Linear algebra approach to sub-problems}
With $l(t)$ we denote the parametric form of a line defined by two distinct point $P$ and $Q$,
$l(t) = (Q-P)t + P$. With $w(S_0E_0, S_1E_1)$ we denote a wedge with arms $S_0E_0$ and $S_1E_1$.
The wedge is called degenerate if its arms are parallel. 
Let 2D cross product of two vectors be defined
as the signed area that these vectors span.

\subsection{Degenerate wedge and line}
\begin{wrapfigure}{r}{0.3\textwidth}
	\vspace{-20pt}
	\centering
	\begin{tikzpicture}[scale=0.4]
	\coordinate [label={above left:$A$}] (A) at (0, 16/3);
	\coordinate [label={below right:$B$}] (B) at (4, 8);
	
	\coordinate [label={right:$S_1$}] (Sr) at (4, 4);
	\coordinate [label={right:$E_1$}] (Er) at (4, 10);
	\coordinate [label={left:$S_0$}] (Sl) at (0, 2);
	\coordinate [label={left:$E_0$}] (El) at (0, 7);
	
	\def\r {2};
	\foreach \o in { (2, 4.26296581635734), (2, 9.070367516975992) } {
		\draw \o circle (\r);	
	}	
	
	\draw[very thick] (Sl)--(El);
	\def\tr{1.3};
     \draw (Er) -- ($ \tr*(Er) - \tr*(Sr) + (Sr)$);
     \def\tr{-0.3};
     \draw (Er) -- ($ \tr*(Er) - \tr*(Sr) + (Sr)$);
     
     \draw [densely dashed, very thick] (A) -- (B);
     
	\draw[very thick] (Sr)--(Er);
	\def\tl{1.9};
     \draw (El) -- ($ \tl*(El) - \tl*(Sl) + (Sl)$);
	
     \def\ra{1.5};
     \draw [densely dashed] (B) -- ($ \ra*(B) - \ra*(A) + (A)$);
     \def\la{1.5};
     \draw [densely dashed] (A) -- ($ \la*(A) - \la*(B) + (B)$);            
\end{tikzpicture}
	\vspace{-10pt}
	\caption{Bootstrapping step}
	\label{fig:DWedgeLine}
	\vspace{-10pt}
\end{wrapfigure}
Let a degenerate wedge $w(S_0E_0, S_1E_1)$ and a line defined by points $A$ and $B$
be given. This case is depicted in Figure~\ref{fig:DWedgeLine}.
Assume that the line is not parallel to the arms and the half of the distance
between the arms equals to $r$. Obviously, inscribed circles will have their centres $I_c$
lying on a line parallel to the arms, such line will be parametrized as, for example,
$l(t)=(E_1-S_1)t+I$, where $I$ is any point equally distanced from the arms.
Radii of the circles will equal to $r$. To find the centres $I_c$ we consider 
the line-to-point distance equation:
\begin{equation}
\frac{(B-A)\times(A-I_c)}{\left\|B-A\right\|}=\pm r.
\end{equation}
After having plugged in the parametric representation of $I_c$, we solve it for $t$
\begin{equation}
t_{1, 2}=\frac{(B-A)\times(A-I)\pm \left\|B-A\right\| r}{(B-A)\times(E_1-S_1)}
\end{equation}
\subsection{Degenerate wedge and point}
\begin{wrapfigure}{r}{0.3\textwidth}
	\vspace{-25pt}
	\centering
	\begin{tikzpicture}[scale=0.4]
	\coordinate [label={right:$S_1$}] (Sr) at (4, 4);
	\coordinate [label={right:$E_1$}] (Er) at (4, 10);
	\coordinate [label={left:$S_0$}] (Sl) at (0, 2);
	\coordinate [label={left:$E_0$}] (El) at (0, 7);
	
	\coordinate [label={above:$P$}] (p) at (2.75, 6.1);
	\fill (p) circle [radius=3pt];
	
	\def\r {2};
	\foreach \o in { (2, 7.954049621773915), (2,  4.245950378226084) } {
		\draw \o circle (\r);	
	}	
	
	\draw[very thick] (Sl)--(El);
	\def\tr{1.1};
     \draw (Er) -- ($ \tr*(Er) - \tr*(Sr) + (Sr)$);
     \def\tr{-0.3};
     \draw (Er) -- ($ \tr*(Er) - \tr*(Sr) + (Sr)$);

	\draw[very thick] (Sr)--(Er);
	\def\tl{1.7};
     \draw (El) -- ($ \tl*(El) - \tl*(Sl) + (Sl)$);
\end{tikzpicture}
	\vspace{-10pt}
	\caption{Bootstrapping step}
	\label{fig:DWedgePoint}
	\vspace{-10pt}
\end{wrapfigure}
Let a degenerate wedge $w(S_0E_0, S_1E_1)$ and a point $P$ between its arms be given.
This case is depicted in Figure~\ref{fig:DWedgePoint}. As before, assume that
the half of the distance between the arms is $r$. Similarly, inscribed circles
will have their centres lying on a line parallel to the arms --- $l(t)=(E_1-S_1)t+I$,
where $I$ is any point equally distanced from the arms. Radii of the circles will equal to $r$.
The centres $I_c$ satisfy the point-to-point distance equation:
\begin{equation}
\left\|P-I_c\right\|^2=r^2
\end{equation}
Substituting $I_c$ with the parametric representation of $l(t)$ brings us to
\begin{equation}
\left\|E_1-S_1\right\|^2 t^2+2(I-P)\cdot(E_1-S_1) t + \left\|I-P\right\|^2-r^2=0
\end{equation}
It is not difficult to see that this equation always has two roots corresponding
to different centres of the inscribed circles.

\eject
\subsection{Wedge and line}
\begin{wrapfigure}{r}{0.4\textwidth}
	\vspace{-25pt}	
	\centering
	\begin{tikzpicture}[scale=0.8]
	\coordinate [label={right:$A$}] (A) at (0.909091, -2.09091);
	\coordinate [label={left:$B$}] (B) at (0., -3.);
	\coordinate [label={above:$C$}] (C) at (0., -1.57143);
	
	\def\r {1.09137};
	\foreach \o in { (1.09137, -3.45206), (2.69118, -1.85226), (-0.872999, -2.32956), (0.726808, -0.729756) } {
		\draw \o circle (\r);	
	}	
	\draw[very thick] (1.09137, -3.45206) circle (\r);
	
     \draw [very thick] (A) -- node [above] {$b$} (C) -- node[left]{$a$} (B);     
     \def\ra{2};
     \draw (B) -- ($ \ra*(B) - \ra*(C) + (C)$);
     \def\la{3};
     \draw (A) -- ($ \la*(A) - \la*(C) + (C)$);
     
     \draw [very thick, densely dashed] (A) -- node[above]{$c$}(B);
     \def\tA{2.2};
     \draw [densely dashed] (B) -- ($ \tA*(B) - \tA*(A) + (A)$);
     \def\tB{2.7};
     \draw [densely dashed] (A) -- ($ \tB*(A) - \tB*(B) + (B)$);
	
	
%
\end{tikzpicture}
	\vspace{-10pt}
	\caption{Typical example where case distinction for~\citep{medvedeva2003implementation}
	 	is required}
	\label{fig:WedgeLine}
	\vspace{-20pt}
\end{wrapfigure}
In Figure~\ref{fig:WedgeLine} a wedge $w(CA,CB)$ and
a line $l$ (dashed) forming $\triangle ABC$ are depicted. We look for an excircle with centre
$I=(x_I, y_I)$ and radius $r$ for $\triangle ABC$ tangent to $AB$.
Let $s$ be a semiperimeter of $\triangle ABC$,
then $r=\sqrt{\frac{s(s-a)(s-b)}{s-c}}$. Let $(x_{AB}, y_{AB})=B-A$ and
$(x_{AC}, y_{AC}) = C-A$. The distance equations between $AB$ and $I$
and $AC$ and $I$ read
\begin{align} 
\frac{(B-A)\times(A-I)}{c} &=  \pm r \\ 
\frac{(C-A)\times(C-I)}{a} &=  \pm r
\end{align}
leading to the following systems of linear equations
\begin{equation}
\begin{pmatrix} 
	-y_{AB} & x_{AB} \\ 
	-y_{AC} & x_{AC}
\end{pmatrix}
\begin{pmatrix} 
	x_I \\ 
	y_I
\end{pmatrix}
=
\begin{pmatrix} 
	B\times A \pm r c\\ 
	C\times A \pm r a\\ 
\end{pmatrix}
\end{equation}
Solution of this system always exists, unless $l$ is parallel to one of the arms of the wedge.
The correct $I$ has to satisfy
\begin{equation}
\left\|\; |(C-B)\times(B-I)|-r b\;\right\|=0.
\end{equation}
In Figure~\ref{fig:WedgeLine} all solutions are plotted and the right excircle
is given in thick.
\subsection{Wedge and point}
\begin{wrapfigure}{r}{0.3\textwidth}
	\vspace{-25pt}	
	\centering
	\begin{tikzpicture}[scale=1.5]
	\coordinate [label={above right:$A$}] (A) at (1.7, 1.8);
	\coordinate [label={left:$B$}] (B) at (1, 1.6);
	\coordinate [label={above:$C$}] (C) at (1, 2.5);
	
	\coordinate [label={below:$P$}] (p) at (1.5, 1.6);
	
	\coordinate (O1) at (1.28998, 1.79994);
	\def\r {0.289975};
	\def\R {0.627182};
	\coordinate (O2) at (1.62718, 0.985848);
		
	\draw [very thick] (A) -- node [above] {$b$} (C) -- node[left]{$a$} (B);
	
	\draw (O1) circle (\r);
	\draw (O2) circle (\R);
	\fill (p) circle [radius=1pt];
	
	\def\ra {2.5};
	\draw [dashed] (B) -- ($ \ra*(B) - \ra*(C) + (C)$);
	
	\def\la {2.5};
	\draw [dashed] (A) -- ($ \la*(A) - \la*(C) + (C)$);
\end{tikzpicture}
	\vspace{-10pt}	
	\caption{Typical example where case distinction for~\citep{medvedeva2003implementation}
		is required}
	\label{fig:WedgePoint}
\end{wrapfigure}
It is known that given a wedge $w(CA, CB)$ and a point $P$ within it,
centres of inscribed circles will lie on bisector of the wedge. The bisector is
defined with points $C$ and $D=(A-B)\frac{a}{a+b}+B$ and thus parametrized
with $l(t)=(D-C)t+C$.
To find the centres $I$ of the inscribed circles we solve
\begin{equation}
\left\|I-P\right\|^2=\left\|\frac{(A-C)\times(C-I)}{b}\right\|^2
\end{equation}
for $t$. Plugging in the parametric representation of the bisector produces
a quadratic equation w.r.t. $t$
\begin{equation}
\left(\left\|D-C\right\|^2-\left(\frac{(A-C)\times(D-C)}{b}\right)^2\right) t^2+
2(D-C)\cdot(C-P) t+\left\|C-P\right\|^2 = 0
\end{equation}
Simple analysis shows, that this equation always has two roots
which correspond to the centres of a smaller and a larger inscribed circles, see 
Figure~\ref{fig:WedgePoint}. According to~\citep{bhattacharya2002minimum} we
are interested in the larger circle.

\section{Conclusion}
The proposed line representation alleviates the fitting procedures removing
the separate treatment for border cases. The parametric representation of lines
guarantees existence of solutions of subsidiary problems.
The algorithm was implemented in TypeScript~\citep{ts} for use
in browsers after transpilation to JavaScript and available on GitHub~\citep{github-release}.


\bibliographystyle{plainnat}      
\bibliography{references}   

\end{document}